\documentclass{beamer}
\usepackage{pstricks}
\newcommand{\dif}{\, \mathrm{d}}

\mode<presentation>
{
  \usetheme{Warsaw}
  \setbeamercovered{transparent}
}

\begin{document}

\begin{frame}
\title
{Multi-scale Analysis  for  Rosseland  Equation\\ with
Small Periodic Oscillating Coef\mbox{}f\mbox{}icients}

\author
{Zhang Qiao-fu}

\institute[]
{ Academy of Mathematics and Systems Science, \\
 Chinese Academy of Sciences  }

\date{NNSFC(No. 90916027)\quad May, \,\,2012}

\maketitle
\end{frame}

 \begin{frame}
\begin{figure}[htbp]
\centering
\includegraphics[height=8cm]{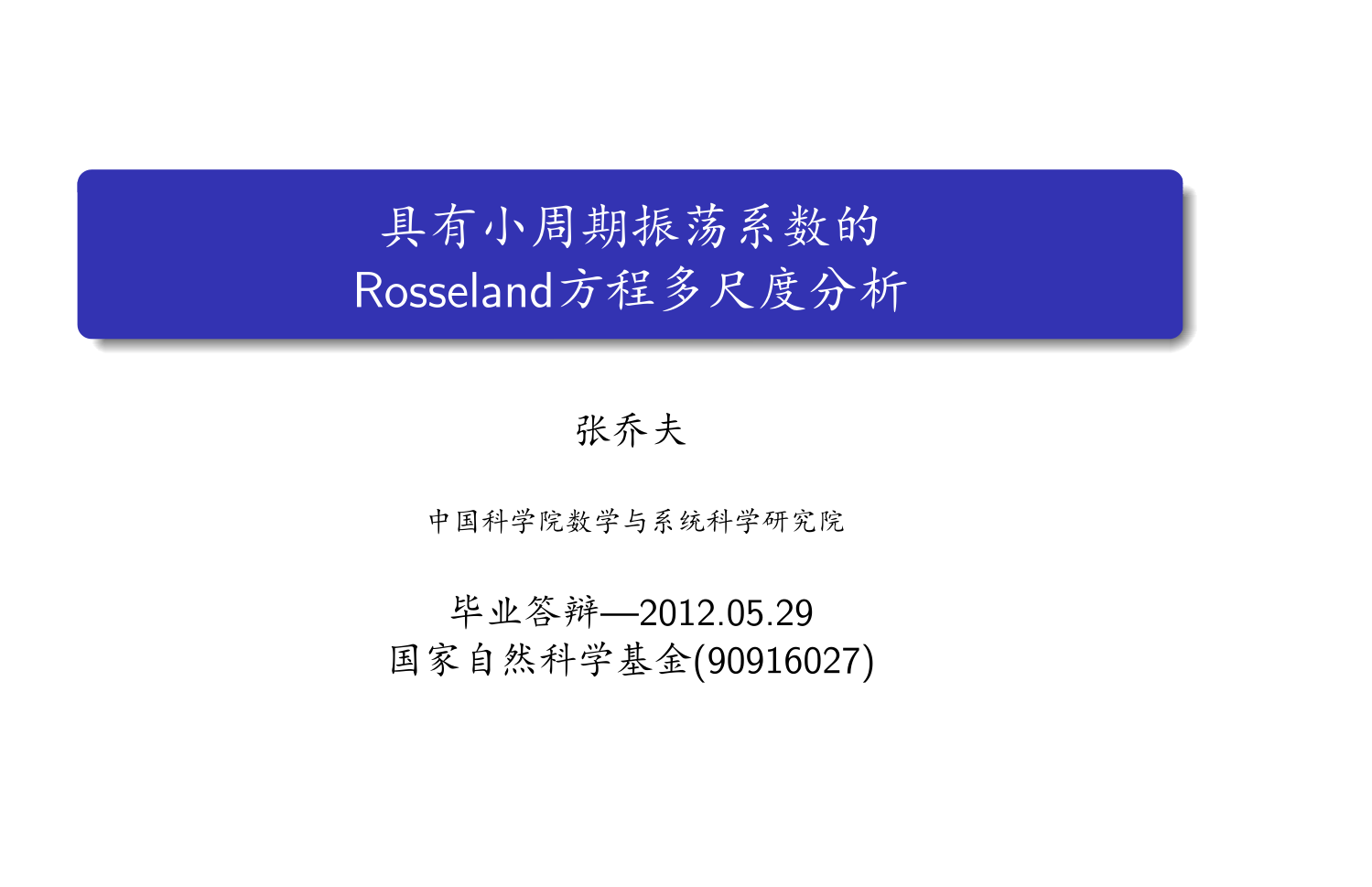}
\end{figure}
\end{frame}


\begin{frame}
  \frametitle{Outline}
 \tableofcontents[pausesections]
\end{frame}

\begin{frame}
  \frametitle{Background (1)}

\begin{figure}[htbp]
\centering
\includegraphics[height=4cm]{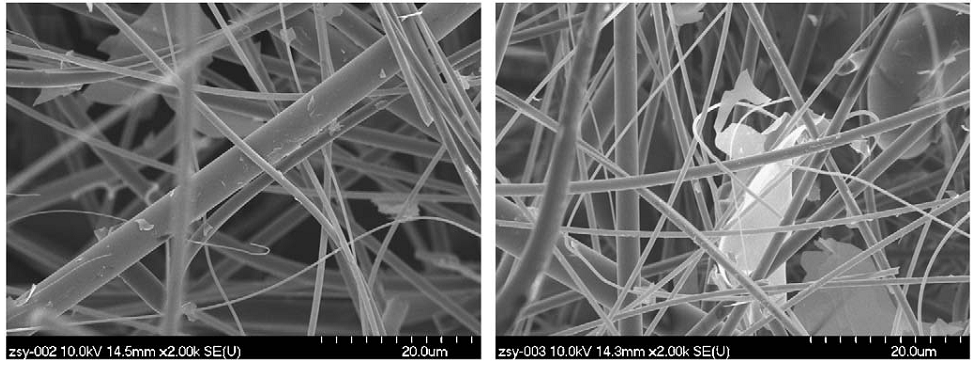}
\end{figure}

  \qquad  SEM of fibrous \quad (H.I.T. 2009) \\
 \pause  \qquad Nonsmooth data

 \end{frame}
 \begin{frame}
  \frametitle{Background (2)}
 \begin{figure}[htbp]
\centering
\includegraphics[height=3cm]{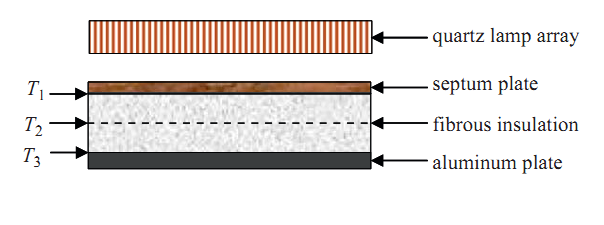}
\end{figure}
   \qquad  \qquad \qquad \qquad  Coat of a rocket \quad  (L.I.R.R., 2005).
\vspace{5mm}

   \end{frame}


 \begin{frame}
  \frametitle{Background (3)}
\begin{figure}[htbp]
\centering
\includegraphics[height=5cm]{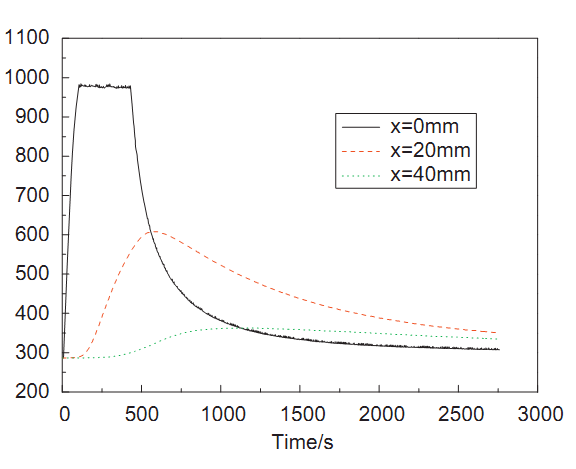}
\end{figure}
  \qquad  \qquad Temperature \,\,(K)\quad  (L.I.R.R., 2005).
\end{frame}

 \begin{frame}
  \frametitle{Background (4)}
\begin{figure}[htbp]
\centering
\includegraphics[height=5cm]{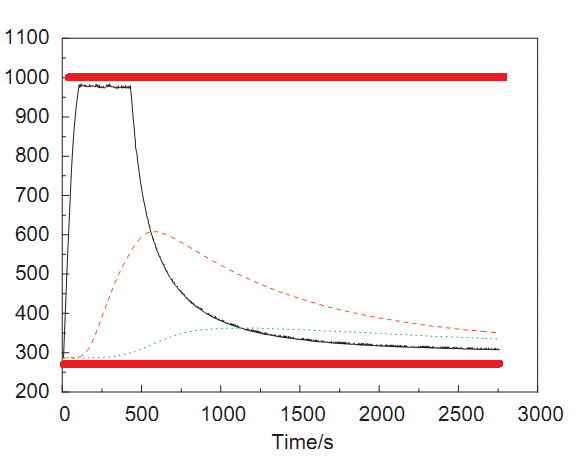}
\end{figure}
 \pause  \qquad  \qquad   No blow-up,\quad  maximum principles,  \pause \quad \textcolor[rgb]{1.00,0.00,0.00}{safe !}

\end{frame}
\section{Rosseland  Eq.---f\mbox{}ixed point}
\subsection{Model}
\begin{frame}
  \frametitle{Part 1: existence / elliptic}

\pause \qquad  \qquad  \,\, Optically \textcolor[rgb]{1.00,0.00,0.00}{thick},   \quad  Rosseland  Eq. \pause
\begin{equation}
 -\mbox{div}[A \nabla\textcolor[rgb]{1.00,0.00,0.00}{u }]=
 f .
\end{equation}
 \qquad  \qquad $u $:\,\, temperature. \pause  \qquad $K$,\,\, $B $ are S P D.
\begin{equation}
A(u ,x ) =K(x )+4  \,\textcolor[rgb]{1.00,0.00,0.00}{u ^3}\,\,  B(x ).
\end{equation}
 \pause  \qquad \qquad Conductive,  \quad radiative.

\end{frame}

 \begin{frame}
  \frametitle{Exists ? }
  \qquad \qquad  Laitinen, 2002
 \begin{figure}
\includegraphics[height=2cm]{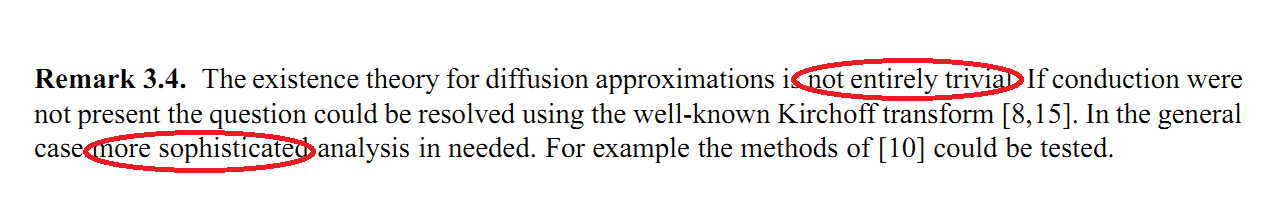}
 \end{figure}
  \pause
 \begin{equation}
A(z,x ) =K(x )+4  \,\textcolor[rgb]{1.00,0.00,0.00}{z^3}\,\,  B(x  ).
\end{equation}
\end{frame}

\begin{frame}
\frametitle{Growth conditions }

\begin{figure}
\includegraphics[height=7cm]{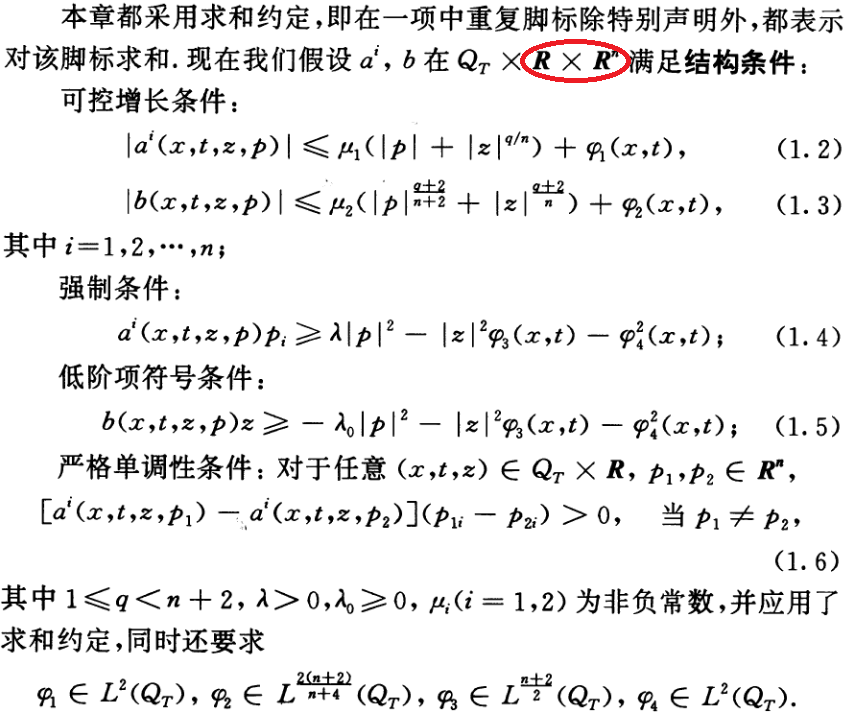}
 \end{figure}
\end{frame}

\begin{frame}
\frametitle{Growth conditions }

\begin{figure}
\includegraphics[height=7cm]{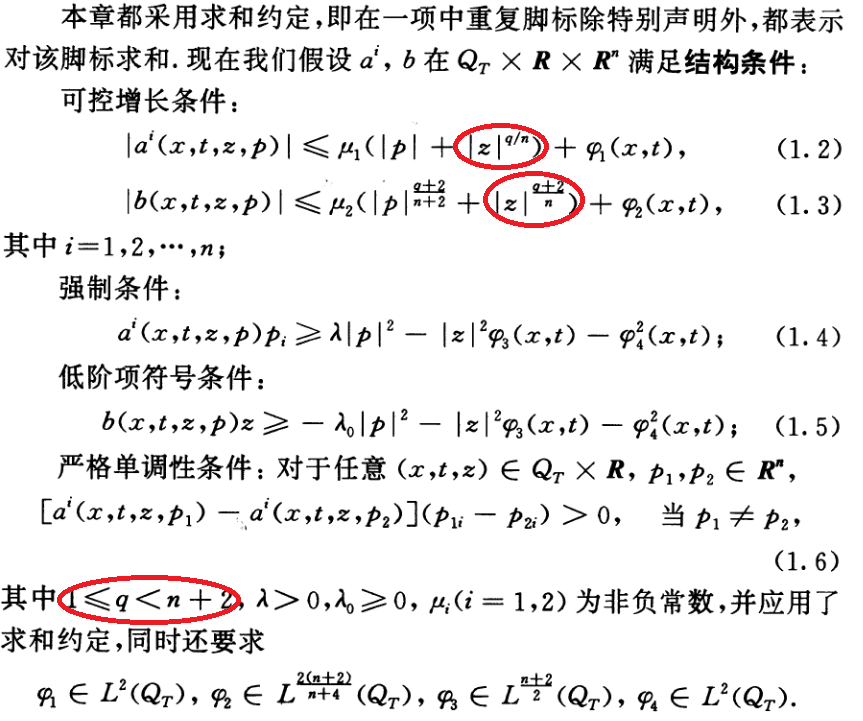}
 \end{figure}

\end{frame}

\begin{frame}
  \frametitle{Physical conditions}
 \pause  \qquad \qquad $0<T_{min}\leq \textcolor[rgb]{1.00,0.00,0.00}{u_b}\leq T_{max}$,
  \pause \quad $0\leq f\leq C$. \pause
\begin{equation}
A(z,x ) =K(x )+4  \,\textcolor[rgb]{1.00,0.00,0.00}{z^3}\,\,  B(x  ).
\end{equation}

\pause
\qquad \qquad If $z\in [T_{min},T_{max} ]$,
\pause
\begin{equation}
\lambda(T_{min},T_{max})\leq A(\textcolor[rgb]{1.00,0.00,0.00}{z} ,x  )\leq \Lambda(T_{min},T_{max} ).
\end{equation}

\end{frame}
\subsection{Elliptic  }
\begin{frame}
  $ \Omega\subset \mathbb{R}^n$,   \quad $u_b \in H^1\cap [T_{min},T_{max} ]$,
\vspace{4mm}

 \pause  F\mbox{}ind $u\in H^{1}(\Omega)$, \quad $(u - u_b )|_{\partial\Omega}=0$, \quad s. t.
\pause \begin{equation}
\int_\Omega A(\textcolor[rgb]{1.00,0.00,0.00}{u(x)},x)\nabla \textcolor[rgb]{1.00,0.00,0.00}{u } \cdot\nabla\varphi = 0,\quad \forall\,\, \varphi\in H^{1}_{0}(\Omega).
\end{equation}

\pause Linearized. \pause \quad F\mbox{}ind $w $, \quad $(w - u_b )|_{\partial\Omega}=0$,  \quad s. t.
\begin{equation}
\int_\Omega A(\textcolor[rgb]{1.00,0.00,0.00}{z(x)},x)\nabla \textcolor[rgb]{1.00,0.00,0.00}{w}  \cdot\nabla\varphi = 0,\quad \forall\,\, \varphi\in H^{1}_{0}(\Omega).
\end{equation}
 \pause
\begin{itemize}
 \pause \item (1)  $\mathcal{L}z=w  \in [T_{min}, \,\,T_{max}]$. \pause \quad $z  \in [T_{min}, \,\,T_{max}] $.
 \vspace{4mm}

 \pause \item (2)  $z=w$?
\end{itemize}

\end{frame}

\begin{frame}
  \frametitle{Linearized map}
\begin{equation}
 \mathfrak{C} =\{ z\in L^2(\Omega);\,\,T_{min}\leq z(x) \leq T_{max},\,\,a.\,e.\,\, in \,\,\Omega\}.
\end{equation}
\pause\qquad    (1)   $\mathfrak{C} $ is closed and convex in $L^2(\Omega)$.

\pause \qquad  (2) $z\in \mathfrak{C} $, \quad $\lambda\leq  A(z)\leq \Lambda$.\pause
\begin{equation}
\int_\Omega A(\textcolor[rgb]{1.00,0.00,0.00}{z(x)},x)\nabla \textcolor[rgb]{1.00,0.00,0.00}{w}  \cdot\nabla\varphi = 0,\quad \forall\,\, \varphi\in H^{1}_{0}(\Omega).
\end{equation}
  \qquad \qquad    $\mathcal{L}z=w $. \pause  \quad    $T_{min}\leq w \leq T_{max}$,\pause \quad $\mathcal{L}\mathfrak{C}\subset \mathfrak{C} $.
\vspace{2mm}

\pause  \qquad  (3) $\|w\|_{H^1}\leq C$, \pause \quad $H^1\Subset L^2$, \quad   $\mathcal{L}\mathfrak{C}  $ is precompact.
\vspace{3mm}

\pause \qquad  (4)  w---\,\,A(z)\pause\,\,---z. \pause \qquad  $\mathcal{L}  $ is continuous.
 \vspace{3mm}

\pause \qquad  \quad A fixed point.
\end{frame}

\begin{frame}
  \frametitle{$A(z)$ is continuous. }
\pause \qquad  $\|z_i-z\|_2\to 0$,\quad $z_i,\,z\in \mathfrak{C}\subset L^\infty$ , \pause \quad $A(z)=K+\textcolor[rgb]{1.00,0.00,0.00}{z^3}\,B$,
\begin{eqnarray}
\|A(z_i)-A(z)\|_2&=&\| (z_i-z)(z_i^2+z_iz+z^2)\,B\|_2
 \nonumber\\
 &\leq&  C \|z_i-z\|_2\to 0 .
 \end{eqnarray}
 \pause
 \begin{eqnarray}
 && \int_{\Omega}|a_{pq}(\textcolor[rgb]{1.00,0.00,0.00}{z_{i }(x)},x)-a_{pq}(\textcolor[rgb]{1.00,0.00,0.00}{z(x)},x)|^2
 \nonumber\\
 &\leq& \int_{\Omega} C |z_{i }(x)-z(x)  |^{2\alpha}
  \nonumber\\
&\leq& C (\int_{\Omega}  |z_{i }(x)-z(x)  |^{\frac{2\alpha}\alpha })^\alpha(\int_{\Omega} 1^{s'})^{\frac 1 {s'}}
  \nonumber\\
&=& C (\int_{\Omega}  |z_{i }(x)-z(x)  |^{\frac{2\alpha}\alpha })^{\frac 1 2 2\alpha}
  \nonumber\\
&\leq&
C\|z_i-z\|_2^{2\alpha}\to 0.
 \end{eqnarray}
\end{frame}

\begin{frame}
  \frametitle{$\mathcal{L}  $ is continuous in $L^2$}
  \pause Suppose $z_i,\,z\in \mathfrak{C},\quad \|z_i-z\|_2\to 0$,\quad\pause $\mathcal{ L}z_i=w_i , \quad \mathcal{ L}z =w.$\pause
  \begin{equation}
\int_\Omega A(\textcolor[rgb]{1.00,0.00,0.00}{z_i})\nabla \textcolor[rgb]{1.00,0.00,0.00}{w_i}  \cdot\nabla\varphi  = 0 ,\quad
\int_\Omega A(\textcolor[rgb]{1.00,0.00,0.00}{z })\nabla \textcolor[rgb]{1.00,0.00,0.00}{w }  \cdot\nabla\varphi  =  0 .
\end{equation}\pause
\begin{equation}
\|z_i-z\|_2\to 0,\quad \|A(\textcolor[rgb]{1.00,0.00,0.00}{z_i})-A(\textcolor[rgb]{1.00,0.00,0.00}{z })\|_2\to 0 .
\end{equation}
 $\|w_i\|_{H^1}\leq C$. \pause \quad $\exists$ $w_0\in H^1( \Omega)$,\quad $(w_0-u_b)\in H^1_0( \Omega)$,
\begin{equation}
 \nabla w_{i_k} \rightharpoonup  \nabla w_0,\quad \mbox{weakly\,\,in\,\,} L^2 ( \Omega; \,\mathbb{R}^n),
\end{equation}
\begin{equation}
  w_{i_k} \to   w_0,\quad \mbox{strongly\,\,in\,\,} L^2 ( \Omega ).
\end{equation}
\quad \quad $\int_\Omega A(\textcolor[rgb]{1.00,0.00,0.00}{z })\nabla \textcolor[rgb]{1.00,0.00,0.00}{w_0}  \cdot\nabla\varphi  =  0
$,  \quad \quad
 $
\int_\Omega A(\textcolor[rgb]{1.00,0.00,0.00}{z })\nabla \textcolor[rgb]{1.00,0.00,0.00}{w }  \cdot\nabla\varphi =  0 .
$\vspace{2mm}

\quad \quad $w_0=w$, \quad $\|  w_{i_k}  -   w_0\|_2\to 0$,  \quad $\|  w_{i }  -   w  \|_2\to 0 $.
\vspace{2mm}

\pause
\quad \quad If $\|z_i-z\|_2\to 0$, \quad  $\|  w_{i }  -   w  \|_2\to 0 $.
\end{frame}

\begin{frame}
  \frametitle{$\mathcal{L}  $ is continuous  in $C^0$}
   \begin{equation}
\int_\Omega A(\textcolor[rgb]{1.00,0.00,0.00}{z_i})\nabla \textcolor[rgb]{1.00,0.00,0.00}{w_i}  \cdot\nabla\varphi  = 0 ,\quad
\int_\Omega A(\textcolor[rgb]{1.00,0.00,0.00}{z })\nabla \textcolor[rgb]{1.00,0.00,0.00}{w }  \cdot\nabla\varphi  =  0 .
\end{equation}\pause \qquad Suppose
\begin{equation}
\|z_i-z\|_\infty \to 0,\quad \|A(\textcolor[rgb]{1.00,0.00,0.00}{z_i})-A(\textcolor[rgb]{1.00,0.00,0.00}{z })\|_\infty\to 0,
\end{equation}
 \pause \qquad then
\begin{equation}
\|w_i-w\|_{C^0(\overline{\Omega})} \to 0,\quad  \|w_i-w\|_{H^1}\to 0.
\end{equation}
 \pause \qquad

 \pause \qquad  Griepentrog,\quad  Recke (2001-2010)
\begin{equation}
\|w_i-w\|_{C^{0,\alpha}(\overline{\Omega})} \to 0,\quad  \|w_i-w\|_{W^{1,2,\omega}}\to 0.
\end{equation}
\pause \qquad Mixed, \quad parabolic systems, \quad nonlinear.
\end{frame}

\begin{frame}
  \frametitle{F\mbox{}ixed point}
  \begin{equation}
\int_\Omega A(\textcolor[rgb]{1.00,0.00,0.00}{z(x)},x)\nabla \textcolor[rgb]{1.00,0.00,0.00}{w}  \cdot\nabla\varphi = 0,\quad \forall\,\, \varphi\in H^{1}_{0}(\Omega).
\end{equation}
 \pause\qquad   (1)  $z\in \mathfrak{C} $, \pause \quad$\mathfrak{C} =\{L^2; [T_{min},T_{max}]\}$
 \pause \quad  \textcolor[rgb]{1.00,0.00,0.00}{closed}  \quad \textcolor[rgb]{1.00,0.00,0.00}{convex.}
\vspace{2mm}

 \pause \qquad   (2) Maximum principles,   \qquad $ \mathcal{L}\mathfrak{C}\subset \mathfrak{C} $.
\vspace{2mm}

\pause \qquad     (3)  $H^1\Subset L^2$,  \qquad $\mathcal{L}\mathfrak{C}  $ is \textcolor[rgb]{1.00,0.00,0.00}{precompact}.
 \vspace{2mm}

\pause \qquad     (4) The solution---data map,   \qquad $\mathcal{L}  $ is \textcolor[rgb]{1.00,0.00,0.00}{continuous}.
 \vspace{2mm}

\pause \qquad  A f\mbox{}ixed point:  $\mathcal{L}u=u  $, \pause \qquad $u$ is a solution.
\vspace{2mm}

\pause \qquad\qquad  $u\in \mathfrak{C}\cap H^1(\Omega)\cap C^{0,\alpha}(\overline{\Omega})$.\pause
\begin{equation}
A\nabla u\cdot\overrightarrow{ n}= \alpha(u-u_{gas}) ,\quad  \alpha\in [0,C].
\end{equation}
\pause \qquad $u_{gas}\in [T_{min},T_{max}]$,  \pause \quad $ \mathcal{L}\mathfrak{C}\subset \mathfrak{C} $.
\end{frame}

\begin{frame}
  \frametitle{Heat source}
   \pause \qquad $0\leq f\leq C $.
 \pause \quad  $\forall\,\textcolor[rgb]{1.00,0.00,0.00}{C_*}>0$,  \pause \quad if $ z \in  [T_{min},\quad  T_{max}+\textcolor[rgb]{1.00,0.00,0.00}{C_*}\,\,]  $, \pause
 \begin{equation}
 0< \lambda (T_{min},T_{max} )
 \leq A(z )\leq  \Lambda (T_{min},T_{max},\textcolor[rgb]{1.00,0.00,0.00}{C_*}\,\,),
\end{equation}
 \pause \qquad
  $z \,\, \in  \,\,\mathfrak{C}  \,\,= \,\,\{L^2\,\,; \quad [T_{min},\quad T_{max} \,\,+ \,\,\textcolor[rgb]{1.00,0.00,0.00}{C_*}\,\,]\,\,\}$, \pause
  \begin{equation}
\int_\Omega A(\textcolor[rgb]{1.00,0.00,0.00}{z } )\nabla \textcolor[rgb]{1.00,0.00,0.00}{w}  \cdot\nabla\varphi
= \int_\Omega f(\textcolor[rgb]{1.00,0.00,0.00}{z})\varphi ,\quad \forall\,\, \varphi\in H^{1}_{0}(\Omega).
\end{equation}
   \pause
 \begin{equation}
  T_{min} \,\,\leq\,\,w \,\,\leq \,\,T_{max}+C_1(\Omega,\lambda(T_{min},T_{max}),\|f\|_\infty).
\end{equation}
\pause \qquad   Let $C_*=C_1$, \pause \quad  $ \mathcal{L}\mathfrak{C}\subset \mathfrak{C} $.
\vspace{2mm}

 \pause \qquad $\lambda (T_{min},T_{max},\textcolor[rgb]{0.00,0.00,1.00}{C_*})$,
 \quad $0\leq f\leq C(T_{min},T_{max},\textcolor[rgb]{0.00,0.00,1.00}{C_*})$.
 \vspace{2mm}

 \pause \quad    $C_*=C_1(C_*)$. \pause   \quad   Non-existence ?
\end{frame}

\subsection{Parabolic }
 \begin{frame}
  \frametitle{Parabolic  model}
\pause\qquad $S=(0,T)$,  \quad    $Q_T=\Omega\times S$, \pause
\begin{equation}
 \textcolor[rgb]{1.00,0.00,0.00}{\partial_t u} -\mbox{div\,}[ A(u,x, t )\nabla u]  = 0 ,\quad \mbox{in\,}\,Q_T .
\end{equation}
\pause
\begin{equation}
 (u-g)|_{t=0}=0   ,\quad (u-g)|_{\partial\Omega}=0 .
\end{equation}
\pause
 \begin{equation}
  g \in H^{1}( Q_T )\cap  [T_{min},\quad T_{max}].
  \end{equation}
 \pause
 \begin{equation}
   \lambda   \leq A(z,x,t) \leq \Lambda  , \quad  z \in  [T_{min},\quad T_{max} ].
  \end{equation}

\end{frame}

 \begin{frame}
  \frametitle{Weak solution}
\qquad
$\forall\, \varphi\in L^2( \,H^1_0 )$,\pause
\begin{equation}
\langle \partial_t u,\,\varphi\rangle_{ L^2( H^1_0 )} +\iint_{Q_T} A(u )\nabla u \cdot \nabla \varphi = 0   .
\end{equation}
\pause
\begin{equation}
\mbox{B.\,\, C.} \quad (u-g)\in \mathcal{W}=\{v\in L^2( H^1_{\textcolor[rgb]{1.00,0.00,0.00}{0}});\, \partial_t v\in L^2( H^{-1})\},
\end{equation}
\pause
\begin{equation}
\mathcal{ W} \hookrightarrow C([0,T];\,L^2(\Omega) ), \quad \mbox{I.\,\, C.} \quad  (u-g)(x,0)=0.
\end{equation}
 \pause \qquad $\langle \partial_t w,\,\varphi\rangle_{ L^2( H^1_0 )}$, \quad integration by parts? \pause
\begin{eqnarray}
&&\langle \partial_t w,\,\varphi\rangle_{ H^1(Q_T )}
\nonumber\\
&=&[\int_\Omega  w \varphi ]|_{0}^T -\iint_{Q_T}  w\partial_t\varphi.
\end{eqnarray}
\pause
\begin{equation}
 w_h(x,\textcolor[rgb]{1.00,0.00,0.00}{t})\,\,\equiv\,\,\frac 1 h \int_{\textcolor[rgb]{1.00,0.00,0.00}{t}}^{\textcolor[rgb]{1.00,0.00,0.00}{t}+h}w(x,s)\,\mbox{d} s.
 \end{equation}
\end{frame}

\begin{frame}
  \frametitle{Linearized map}
\pause \begin{equation}
 \mathfrak{C} =\{ \phi\in L^2(Q_T);\,\,T_{min}\leq \phi(x,t) \leq T_{max},\,\,\mbox{a.\,e.\,\, in} \,\,Q_T\}.
\end{equation}
\pause\qquad  is a closed convex  set in   $L^2(Q_T )$.
\vspace{2mm}

\pause \qquad $\forall\, \varphi\in L^2(S;\,H^1_0(\Omega))$, \pause
\begin{equation}
\langle \partial_t w,\,\varphi\rangle_{ L^2(S;\,H^1_0(\Omega))} +\iint_{Q_T} A( \textcolor[rgb]{1.00,0.00,0.00}{z} )\nabla w \cdot \nabla \varphi = 0   .
\end{equation}
\pause
\begin{equation}
(w-g)\in \mathcal{W} ,\qquad (w-g)(x,0)=0.
\end{equation}

\pause \qquad$\mathcal{ L}z=w$. \pause \qquad $\mathcal{ L} \mathfrak{C} \subset  \mathfrak{C}$?
\end{frame}

\begin{frame}
\frametitle{Maximum principles (1)}
\pause \qquad $\mathcal{W}=\{v\in L^2( H^1_0);\, \partial_t v\in L^2( H^{-1})\}$ \pause\,\,=\,\, $\overline{C^\infty([0,T];\,H^1_0(\Omega))} $.
\vspace{5mm}

\pause $\forall\, (w-g)\in \mathcal{W}$, \pause\quad  $\exists \,\textcolor[rgb]{1.00,0.00,0.00}{\{v_i\}}
\subset C^\infty( H^1_0 )$, \pause\quad $\|w-g-\textcolor[rgb]{1.00,0.00,0.00}{v_i}\|_{\mathcal{W}}\to 0$,
\pause
\begin{eqnarray}
&&\qquad \,\,\,\,\,\,\langle \partial_t (v_i+g),\,(v_i+g-T_{max})_+\rangle_{ L^2( H^1_0 )}
 \nonumber\\
 &&\qquad= \int_0^T\langle \partial_t (v_i+g),\,(v_i+g-T_{max})_+\rangle_{  H^1_0  }
 \nonumber\\
&&\qquad=
 \iint_{Q_T} \partial_t (v_i+g)\cdot(v_i+g-T_{max})_+
\end{eqnarray}
\pause \begin{eqnarray}
&=\,\,&
 \iint_{Q_T} \partial_t  \frac {(v_i+g-T_{max})_+^2(x,t)}{2}
   \nonumber\\
 & =\,\, & \int_\Omega  \frac {(v_i+g-T_{max})_+^2(x,T) }{2}-0
     .
\end{eqnarray}
\end{frame}

 \begin{frame}
  \frametitle{Maximum principles (2)}
\begin{equation}
 \langle \partial_t w,\,(w-T_{max})_+\rangle_{ L^2( H^1_0 )}
=\int_\Omega  \frac {(w-T_{max})_+^2 (x,T) }{2}\geq 0
     .
\end{equation}\pause
\begin{eqnarray}
 &&\iint_{Q_T} A(z )\nabla w \cdot \nabla (w-T_{max})_+
 \nonumber\\
 &=&
  \iint_{Q_T} A(z )\nabla (w-T_{max})_+ \cdot \nabla (w-T_{max})_+
  \nonumber\\
 & \geq &
  \lambda  \int_S\int_{\Omega} |\nabla (w-T_{max})_+|^2
   \nonumber\\
 & \geq &
  C(\Omega)\lambda \int_S\int_{\Omega}   (w-T_{max})_+ ^2
   .
\end{eqnarray}
\end{frame}

 \begin{frame}
  \frametitle{Maximum principles (3)}
\pause \qquad Let $\varphi=(w-T_{max})_+\in L^2( H^1_0 )$, \pause
\begin{eqnarray}
 0
 & = &  \langle \partial_t w,\,(w-T_{max})_+ \rangle_{ L^2(S;\,H^1_0(\Omega))}
  \nonumber\\
  &&\quad\quad
  +\iint_{Q_T} A(z )\nabla w \cdot \nabla (w-T_{max})_+
   \nonumber\\
 &\geq & 0\,\,+\,\, C(\Omega)\lambda \int_S\int_{\Omega}   (w-T_{max})_+ ^2 \,\,\geq\,\, 0.
\end{eqnarray}
\pause \qquad  $w\leq T_{max}$.  \pause \qquad $w\geq T_{min}$,\quad  $\textcolor[rgb]{1.00,0.00,0.00}{\mathcal{L} \mathfrak{C }\subset\mathfrak{C }}$.
\pause
\begin{equation}
\|w\|_{\mathcal{W }}^2=\|w\|_{ L^2( H^1_0 )}^2+\|
 \partial_t w \|_{ L^2( H^{-1})}^2
 \leq C.
\end{equation}
\pause\qquad $\mathcal{W } \Subset L^2(Q_T)$, \pause \qquad
$\mathcal{L} \mathfrak{C }$ is \textcolor[rgb]{1.00,0.00,0.00}{precompact}.
\vspace{2mm}

\pause\qquad A fixed point if $\mathcal{L}$ is \textcolor[rgb]{1.00,0.00,0.00}{continuous}.
\end{frame}

 \begin{frame}
  \frametitle{ $\mathcal{L}$ is continuous.}
\pause
\qquad
$\forall\, \varphi\in L^2( H^1_0 )$, \quad $z_i,\,z\in\mathfrak{C }$,\pause
\begin{eqnarray}
\langle \partial_t w_i,\,\varphi\rangle_{ L^2( H^1_0 )} +\iint_{Q_T} A(\textcolor[rgb]{1.00,0.00,0.00}{z_i})\nabla w_i \cdot \nabla \varphi& =& 0   ,
\\
\langle \partial_t w ,\,\varphi\rangle_{ L^2( H^1_0 )} +\iint_{Q_T} A(\textcolor[rgb]{1.00,0.00,0.00}{z} )\nabla w  \cdot \nabla \varphi &=& 0   .
\end{eqnarray}
\pause
\qquad  $\{w_i\}$ is bounded in $\mathcal{W }$. \pause
\qquad $\partial_t w_i\rightharpoonup  \partial_t w_0$, \quad  weakly in $L^2(H^{-1}) $.
\vspace{2mm}

\qquad $\|A(z_i)-A(z)\|_2\to 0$,
\qquad $\nabla w_i \rightharpoonup \nabla w_0$,\quad  weakly in $L^2$. \pause
\begin{equation}
\langle \partial_t w_0,\,\varphi\rangle_{ L^2( H^1_0 )} +\iint_{Q_T} A(z )\nabla w_0 \cdot \nabla \varphi = 0   .
\end{equation}
\qquad $w_0=w$.
 \qquad $\mathcal{W } \Subset L^2(Q_T)$,\pause  \quad $\|w_i-w\|_2\to 0$.
\end{frame}

 \begin{frame}
  \frametitle{Integral }
\pause \qquad If $u|_{[-h_0,0]}$ is known, \pause
\begin{equation}
  \partial_t w-\mbox{div}[ A(\frac 1 {h_0}\int^t_{t-h_0}\textcolor[rgb]{1.00,0.00,0.00}{z(x,s)}\dif s  )\nabla w ]  = 0   .
\end{equation}
\pause
\begin{equation}
\langle \partial_t w ,\,\varphi\rangle_{ L^2( H^1_0 )} +\iint_{Q_T} A(\frac 1 {h_0}\int^t_{t-h_0}\textcolor[rgb]{1.00,0.00,0.00}{z(s)}\dif s  )\nabla w  \cdot \nabla \varphi = 0   .
\end{equation}
\pause \qquad Similarly,

\end{frame}

 \begin{frame}
  \frametitle{Elastic wave ...}
  \begin{figure}
\includegraphics[height=2cm]{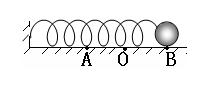}
 \end{figure}
\pause  \qquad $u_{max}$ is independent of $ \lambda$, $ \Lambda$. \pause
\begin{equation}
\langle  w'' ,\,\varphi\rangle_{ L^2( H^1_0 )} +\iint_{Q_T} A_{ij}^{hk}( z ,x, t )\frac {\partial w_i}  {\partial x_k}  \frac {\partial  \varphi_j}  {\partial x_h}= 0   .
\end{equation}
 \pause  \qquad $L^\infty$ ?
\end{frame}

 \begin{frame}
  \frametitle{Wave ...}
\begin{equation}
\langle  w'' ,\,\varphi\rangle_{ L^2( H^1_0 )} +\iint_{Q_T} A(  x, t )\nabla w\cdot \nabla \varphi =
\iint_{Q_T} f(\textcolor[rgb]{1.00,0.00,0.00}{z},x,t) \varphi   .
\end{equation}
\pause
\begin{equation}
\| w\|_{C(H^1_0)} +\|w'\|_{L^\infty(L^2) } + \|w''\|_{L^2(H^{-1})} \leq C_1(\|f\|_2+\|\phi_0\|_{H^1}+\|\phi_1\|_{L^2} ) .
\end{equation}
\pause  \qquad Suppose (1) $f(z,x,t)\in L^2$ for any $z$, \quad  e.g. $f\in L^\infty$.
\vspace{2mm}

\pause \qquad  (2) $n=1$, \quad $ C(H^1_0) \hookrightarrow C^0(\overline{Q}_T)$.
\vspace{2mm}

\pause  \qquad  (3) $A(x,t)\in [ \lambda, \Lambda]$,\quad  $ z\in \mathfrak{C}= \{ L^2(Q_T)\,;\,\, [-C_1,\,\, C_1]\}$.
\pause   \quad $\textcolor[rgb]{1.00,0.00,0.00}{\mathcal{L} \mathfrak{C }\subset\mathfrak{C }}$.
\vspace{2mm}

 \pause  \qquad $H^1(Q_T) \Subset L^2(Q_T)$. \pause \qquad
$\mathcal{L} \mathfrak{C }$ is \textcolor[rgb]{1.00,0.00,0.00}{precompact}.
\vspace{2mm}

\pause\qquad A f\mbox{}ixed point if $\mathcal{L}$ is \textcolor[rgb]{1.00,0.00,0.00}{continuous}.
\end{frame}

 \section{Math algorithms}
\begin{frame}
  \frametitle{Part 2: Math algorithms}
\pause \qquad$Q_t=\Omega\times (0,t)$. \pause\quad If  t is small enough,\pause
\begin{eqnarray}
\langle   \partial_t u_{i+1},\,  \varphi \rangle_{W^{1,1}_D(Q_t)}
&+&\iint_{Q_t}  A(\textcolor[rgb]{1.00,0.00,0.00}{ u_i} )\nabla u_{i+1}\cdot\nabla \varphi \nonumber\\
+\iint_{(0, t)\times \Gamma} \alpha(u_{i+1}-u_{gas})\varphi
&=&\iint_{Q_t}  f(\textcolor[rgb]{1.00,0.00,0.00}{u_i} )\varphi  .
\end{eqnarray}
\pause \qquad $A(u_{i+1}) \approx A(u_{i })$.
\vspace{4mm}

\pause\qquad $A(u_{i+1}) \approx A(u_{i })+A'(u_i)(u_{i+1}-u_{i } )$, \quad Newton method.
\vspace{4mm}

\pause\qquad  $\partial_t u\approx \frac {u(t_{j+1})-u(t_j)} h$,\quad Rothe method.
\end{frame}

 \section{SOTS---homogenization}
\begin{frame}
\frametitle{Part 3: SOTS---homogenization}
\pause
\begin{equation}-\mbox{div}[A(\textcolor[rgb]{1.00,0.00,0.00}{u_\varepsilon} ,x,\frac x \varepsilon)\nabla u_\varepsilon]
=0.
\end{equation}
\pause \begin{equation}
u_{\varepsilon}( x)=u_0(x) +
\varepsilon u_1(x,\frac{x}{\varepsilon}) + {\varepsilon}^2
u_2(x,\frac{x}{\varepsilon}) + ...
\end{equation}
\pause \begin{equation}
A(u_{\varepsilon})=A(u_0,x) + \varepsilon A_1(u_0,u_1,x,\frac{x}{\varepsilon}) + {\varepsilon}^2A_2  + ...
\end{equation}
\pause
\begin{equation}
\int_\Omega A_0(u_0 ,x)\nabla u_0 \cdot\nabla \varphi
=0.
\end{equation}
\pause
\begin{equation}
  a_{ij}^0(u_0,x)
=\int_Y[a_{ij}(u_0 )+a_{il}(u_0 )\frac{\partial }{\partial y_l}N_j(u_0,x,y)]\mbox{d} y.
\end{equation}
\pause
\begin{equation}
\int_Y A(\textcolor[rgb]{1.00,0.00,0.00}{ u_0 })\nabla (N_j +y_j)\cdot\nabla \varphi=0,\quad\forall\, \varphi
\in W^1_{per}(Y) .
\end{equation}
\end{frame}

\begin{frame}
\frametitle{Homogenization (1)}
\qquad If $z\in [T_{min},T_{max} ] $, \pause\qquad $\lambda\leq A(z)\leq \Lambda $.
\vspace{4mm}

\pause\qquad The div-curl Lemma,
 \pause \qquad $\lambda\leq A_0(z)\leq \Lambda_0$.
\vspace{4mm}

\pause\qquad If $A(z)$ is continuous,  \pause \qquad $  A_0(z) $ is also continuous.
\pause
\begin{equation}
\exists u_0,\quad \mbox{s.\,t.}\quad \int_\Omega A_0(u_0 ,x)\nabla u_0 \cdot\nabla \varphi
=0.
\end{equation}
 \pause \qquad $\|u_\varepsilon-u_0\|_{C^0(\overline{\Omega})} \to 0$,
  \vspace{4mm}

  \qquad  $A (u_\varepsilon)\nabla u_\varepsilon  \rightharpoonup A_0(u_0)\nabla u_0$, \quad weakly in $L^2$.
\vspace{4mm}

\pause\qquad Fusco, Moscariello,  \quad $A_{i\nu}\equiv A(\langle w(x)\rangle_{i\nu},\langle  x \rangle_{i\nu},y)$.
\end{frame}

\begin{frame}
\frametitle{Homogenization (2)}
\pause \qquad  $\forall\,\varepsilon>0$, $\exists \, u_\varepsilon $,  s.t.
\begin{equation}
\int_\Omega A (\textcolor[rgb]{1.00,0.00,0.00}{u_\varepsilon} ,\frac{x}{\varepsilon})\nabla u_\varepsilon \cdot\nabla \varphi
=0.\qquad \|u_{\varepsilon_k}-u_0\|_{C^0}\to 0.
\end{equation}
\pause \begin{equation}
\int_\Omega A (\textcolor[rgb]{1.00,0.00,0.00}{u_0} ,\frac{x}{\varepsilon})\nabla u_{0,\varepsilon} \cdot\nabla \varphi
=0,
\qquad \int_\Omega A_0 (\textcolor[rgb]{1.00,0.00,0.00}{u_0 } )\nabla w_0 \cdot\nabla \varphi
=0.
\end{equation}
\pause \begin{equation}
\| u_{0,\varepsilon}-w_0\|_{L^2}\to 0,\quad
A (u_0 ,\frac{x}{\varepsilon})\nabla u_{0,\varepsilon}
\rightharpoonup A_0 (u_0 )\nabla w_0  .
\end{equation}
\pause \qquad$\| u_\varepsilon  -u_0  \|_{C^0}\to 0$,\quad $\| A (u_\varepsilon ,\frac{x}{\varepsilon})- A (u_0 ,\frac{x}{\varepsilon})\|_{L^\infty}\to 0$,
\vspace{3mm}

\pause \quad  $\|u_\varepsilon-u_{0,\varepsilon}\|_{H^1}\to 0   $.\pause \quad  $\|u_\varepsilon-u_{0,\varepsilon}\|_{L^2}\to 0   $. \pause \quad $u_0=w_0$.
\vspace{3mm}


 \pause \qquad$ A (u_\varepsilon ,\frac{x}{\varepsilon})\nabla u_\varepsilon  \approx A (u_0 ,\frac{x}{\varepsilon})\nabla u_{0,\varepsilon}$  \pause $\rightharpoonup   A_0 (u_0 )\nabla u_0  $.
\end{frame}

\section{Error in $C^{0,\alpha}_{loc}$ and $W^{1,\infty}_{loc}$}
\subsection{ Correctors estimates}
\begin{frame}
\frametitle{Part 4: Error in $C^{0,\alpha}_{loc}$ and $W^{1,\infty}_{loc}$}
 \vspace{4mm}
\pause \qquad $Y=(0,1)^n$,\pause\quad $W^1_{per}(Y)=\{H^1(Y):$ Y-periodic, $\int_Y \varphi =0\}$.
 \vspace{4mm}

\pause F\mbox{}ind $P\in W^1_{per}(Y)$, s.t.
\begin{equation}
  \int_{Y} A \nabla P\cdot \nabla\varphi
  =\int_{Y} \overrightarrow{ B} \cdot \nabla\varphi+\int_{Y}d \varphi,\quad\forall\, \varphi\in W^1_{per}(Y).
\end{equation}
\pause \qquad $A$ is S P D.\pause \quad  $A$, $\overrightarrow{ B},\,\,d $ are Y-periodic.  \qquad $\int_Y d =0$.
 \vspace{4mm}

 \pause \qquad Lax-Milgram Lemma in $ W^1_{per}(Y)$.
\end{frame}

\begin{frame}
\frametitle{Translation invariance }
\begin{equation}
  \int_{Y} A \nabla \textcolor[rgb]{1.00,0.00,0.00}{P}\cdot \nabla\varphi
  =\int_{Y} \overrightarrow{ B} \cdot \nabla\varphi+\int_{Y}d \varphi,\quad\forall\, \varphi\in W^1_{per}(Y).
\end{equation}
\pause \qquad Boundary estimates ? \quad Maximum principles ?
\begin{figure}[htbp]
 \centering
\includegraphics[ width= 2cm ]{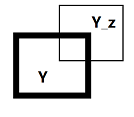}
\end{figure}

\pause \qquad $z\in \mathbb{R}^n$,\quad  $ Y_z=Y+z$,\pause
\begin{equation}
  \int_{Y_z} A \nabla \textcolor[rgb]{1.00,0.00,0.00}{Q}\cdot \nabla\varphi
  =\int_{Y_z} \overrightarrow{ B } \cdot \nabla\varphi+\int_{Y_z}d \varphi,\quad\forall\, \varphi\in W^1_{per}(Y_z).
\end{equation}
\end{frame}

\begin{frame}
\frametitle{Local---Global}
 \qquad If $\zeta \,\,\textcolor[rgb]{1.00,0.00,0.00}{( =A\nabla P- \overrightarrow{ B })}$ \pause $\in L^2_{per}(Y;\mathbb{R}^n)$,
 \pause \quad $d \in L^2_{per}(Y )$, \quad $\int_Y d  =0$,\pause
\begin{equation}
  \int_Y  \zeta(y)\cdot \nabla \varphi + d(y)\varphi =0,\qquad \forall \,\varphi(y)\in W^1_{per}(Y).
\end{equation}
\pause \qquad For any bounded Lipschitz domain $ K\subset \mathbb{R}^n$, \pause
\begin{equation}
\int_{\textcolor[rgb]{1.00,0.00,0.00}{K }}\zeta(x)\cdot \nabla \phi(x)+ d(x) \phi(x) =0,\quad \forall \,\phi(x)\in\textcolor[rgb]{1.00,0.00,0.00}{ H^{1}_{0}(K)}.
\end{equation}
\pause \qquad   If $d=0$,
\begin{equation}
-\mbox{div}\,  _x\zeta(x)=0, \quad \mbox{ in}\,\,   H^{-1}(K).
\end{equation}
\pause \qquad  (1986)$\to$(1982),\quad div-curl. \pause
 \begin{equation}
  \int_{Y} A \nabla P\cdot \nabla\varphi
  =\int_{Y} \overrightarrow{ B} \cdot \nabla\varphi+\int_{Y}d \varphi,\quad\forall\, \varphi\in W^1_{per}(Y).
\end{equation}
\pause \qquad  $\|P\|_{C^{0,\alpha}(\overline{Y})}$. \pause \qquad $\|P\|_{W^{1,\infty}( Y )}$?
\end{frame}

\begin{frame}
\frametitle{Gradient estimates}
 \pause \qquad Li and Vogelius (2000).\quad
\quad Divergence form elliptic equations, \pause\quad piecewise H\"older continuous coef\mbox{}f\mbox{}icients.
\vspace{3mm}

 \pause \qquad   Global $W^{1,\infty}$ and piecewise $C^{1,\alpha}$ estimates.   \pause
\begin{equation}
  \int_{Y} A \nabla P\cdot \nabla\varphi
  =\int_{Y} \overrightarrow{ B} \cdot \nabla\varphi+\int_{Y}d \varphi,\quad\forall\, \varphi\in W^1_{per}(Y).
\end{equation}
 \pause \qquad $\|P\|_{W^{1,\infty}( Y )} \leq C$. \pause \qquad C. D. (1999). \pause
\begin{equation}
 u_0(x)+\varepsilon N_k(\textcolor[rgb]{1.00,0.00,0.00}{\frac x \varepsilon})\frac{\partial   u_0}{\partial x_i} +
 \varepsilon^2 M_{kl}(\textcolor[rgb]{1.00,0.00,0.00}{\frac x \varepsilon}) \frac{\partial^2 u_0}{\partial x_k \partial x_l  }.
\end{equation}
\end{frame}
\subsection{ Error estimates}
\begin{frame}
\frametitle{Error estimates (1)}
\begin{equation}
Z_\varepsilon\equiv u_\varepsilon-u_0(x)-\varepsilon N_k(\textcolor[rgb]{1.00,0.00,0.00}{\frac x \varepsilon})\frac{\partial   u_0}{\partial x_i} -
 \varepsilon^2 M_{kl}(\textcolor[rgb]{1.00,0.00,0.00}{\frac x \varepsilon}) \frac{\partial^2 u_0}{\partial x_k \partial x_l  }.
\end{equation}
\pause
\begin{equation}
   -\frac{\partial}{\partial x_i}\left(a_{ij}( \frac x \varepsilon)\frac{\partial Z_\varepsilon}{\partial x_j}
     \right)=\varepsilon \frac{\partial \theta_i}{\partial x_i} + \varepsilon \eta,
\end{equation}
 \pause \qquad If $u_0\in W^{3,q}$, \quad piecewise H\"older continuous,\pause\quad $ \theta_i$, $\eta \in L^q$.
 \vspace{3mm}

 \pause \qquad  Maximum principles,\quad De Giorgi-Nash estimates, \pause
\begin{equation}
\sup_{\Omega} |u_\varepsilon-u_0|\leq C\varepsilon,\quad \mbox{B.\,\,L.\,\,P.\,\,\,\,1978},
\end{equation} \pause
  \begin{equation}
    \|u_\varepsilon-u_0-\varepsilon u_1\|_{C^{0,\beta}(\overline{\Omega'})}\leq C\varepsilon, \quad \Omega'\subset\subset\Omega;
\end{equation} \pause
\begin{equation}
  |\int_\Omega A(\frac x\varepsilon)\nabla u_\varepsilon\cdot\nabla u_\varepsilon-
  \int_\Omega A_0\nabla u_0\cdot\nabla u_0|\leq C\varepsilon;
\end{equation}
 \pause \qquad Error in $W^{1,\infty}$ ?
\end{frame}

\begin{frame}
\frametitle{Error estimates (2)}
Avellaneda-Lin (1987)
\begin{equation}
   -\mbox{div}[A(\frac x\varepsilon )\nabla u_\varepsilon]
   =f,
\end{equation}
 \pause \qquad $A(y)$ is S P D, Y-periodic and $C^{0,\gamma}(\overline{Y})$ (or p-w smooth),  \quad $ f\in L^q$, \pause
\begin{equation}
   \|\nabla u_\varepsilon\|_\infty\leq C(\|u_\varepsilon\|_\infty+ \|f\|_q ).
\end{equation}
 \pause \qquad If $u_0\in W^{4,q}_{loc}$,  \pause \quad $A(y)\in C^{0,1}(\overline{Y})$,\pause\quad  R.H. of $Z_\varepsilon\in L^q$, \pause
\begin{equation}
\sup_{\overline{\Omega'}}|\nabla (u_\varepsilon-u_0-\varepsilon u_1)|\leq
 C\varepsilon,\quad \Omega'\subset\subset \Omega.
\end{equation}
\end{frame}

\begin{frame}
\frametitle{ Conclusions }
\begin{itemize}
      \item (1) PD in an interval,\quad nonsmooth f\mbox{}ixed point.
\vspace{3mm}

\pause  \item (2) Mathematical algorithms.
\vspace{3mm}

\pause  \item (3) SOTS --- homogenization.
\vspace{3mm}

\pause  \item (4) Error estimates   in $C^{0,\alpha}_{loc}$ and $W^{1,\infty}_{loc}$.
\vspace{5mm}
\end{itemize}
 \pause \quad \quad  (a) Temperature,  \pause  \quad  (b) $C^{0,\alpha}$ and p-w smooth --- $W^{1,\infty}$.
\end{frame}

\begin{frame}
\frametitle{Problems }
\begin{itemize}
  \pause      \item (1)
Gradient estimates,\quad Li-Vogelius (p-w),

\pause \qquad Avellaneda-Lin ( $A(x, \frac x\varepsilon)$ ).
\vspace{3mm}

\pause  \item (2) F\mbox{}lux  estimates, \pause \quad $A\nabla u$,   $A\nabla u_h$,

\pause \qquad  p-w,\quad $\|A _\varepsilon\nabla u_\varepsilon\cdot \overrightarrow{n }-A _0\nabla u_0\cdot \overrightarrow{n }\|_{C^0}\to 0$ ?
\vspace{3mm}

\pause  \item (3)  Uniqueness.
\vspace{3mm}

\pause  \item (4) Nonlinear data--solution,\quad $L^\infty$ estimates for hyperbolic eq.
\vspace{5mm}
\end{itemize}

\vspace{3mm}
\pause\quad\quad\quad\quad \Large{\textcolor[rgb]{1.00,0.00,0.00}{Thanks !}}

\end{frame}

\end{document}